%
\input amstex \documentstyle{amsppt} \magnification=\magstep1
\hsize=6true in \vsize=9true in \NoBlackBoxes \TagsOnRight
\NoRunningHeads \UseAMSsymbols  \nologo   

\def\section#1{\heading{#1}\endheading}
\def\pf{\demo{Proof}} \def\endpf{\quad$\square$\enddemo} 

 \def\={\;=\;} \def\:{\;:=\;}

\def\a{\alpha} \def\b{\beta} \def\C{\Bbb C}  \def\D{\Delta}
   \def\eq{$(\#)_k$ } 
\def\eqm{$(\#)_{k-6}$ }  \def\eqp{$(\#)_{k+6}$ }  \def\G{\Gamma}
\def\g{\text{SL}_2(\Z)}  \def\h{\frak{H}}  \def\i{^{-1}}
\def\inf{\infty}   \def\j1{j(\tau)} \def\j2{j_{(2)}(\tau)}
\def\j3{j_{(3)}(\tau)} \def\j4{j_{(4)}(\tau)} 
\def\l{\lambda}    \def\O{\text{\rm O}}
   \def\p{\partial}  \def\Q{\Bbb Q} 
  \def\s{\sigma}  \def\t{\tau} 
 \def\Z{\Bbb Z}

\topmatter\title On modular forms arising from a differential \\
equation of hypergeometric type\endtitle \author Masanobu Kaneko and
Masao Koike\endauthor \endtopmatter\document

\section{\S1. Introduction}

The differential equation that the present paper shall deal with is
$$(\#)_k\qquad f''(\tau)-\frac{k+1}6
E_2(\tau)f'(\tau)+\frac{k(k+1)}{12} E'_2(\tau)f(\tau)=0,$$
where
$\tau$ is a variable in the upper half-plane $\h$, the symbol ${}'$ a
differential operator   $(2\pi i)\i d/d\tau=q\cdot d/dq\ (q=e^{2\pi i
  \tau})$, and $E_2(\tau)$ the ``quasimodular'' Eisenstein series of
weight $2$ for the full modular group $\g$:
$$E_2(\tau)=1-24\sum_{n=1}^\inf\bigl(\sum_{d|n} d\bigr)q^n.$$
Note
that our usage of the symbol ${}'$ differs from the standard one by
the  factor $1/2\pi i$. The parameter $k$ always stands for a
nonnegative integer or  half an integer throughout the paper.   This
differential equation originates in the work [1] where in some cases
($k\equiv 0, 4 \mod 6$) solutions which are modular on $\g$ were found
and studied in connection with liftings of supersingular
$j$-invariants of elliptic curves.

The purpose of this paper is to give an explicit description  of
(conjecturally) all modular solutions of \eq when $k$ is an integer or
half an integer (Theorem 1 in Section 2), as well as to discuss
positiveness of Fourier coefficients of some of those solutions
(Theorem 3 in Section 4).  We shall also  discuss an intrinsic
characterization of the equation \eq by the property that if $f(\t)$
is a solution, then $(c\t+d)^{-k}f((a\t+b)/(c\t+d))$ is also a
solution for all $\left(\smallmatrix
  a&b\\c&d\endsmallmatrix\right)\in\g$ (Proposition 2 in Section 5).
When $k$ is an odd integer congruent to $5$ modulo $6$, an unexpected
solution occurs in contrast to the other cases:  \eq has {\it
  quasimodular} solutions of weight $k+1$ (rather than $k$ as in other
modular solutions).  We shall describe this quasimodular solution in
terms of certain orthogonal polynomials (Theorem 2 in Section 2), and,
to show that this is a solution, discuss an inductive structure of
solutions of  \eq with different $k$'s in Section 3.
 
\section{\S2. Explicit description of modular and quasimodular solutions}

To describe the solutions, we need to develop notations of  various
modular forms of levels $1,2,3,$ and $4$.  Let
$$E_4(\t)=1+240\sum_{n=1}^\inf \bigl(\sum_{d|n} d^3\bigr)q^n=1 + 240q + 2160q^2 
+ 6720q^3 +\cdots $$ and
$$E_6(\t)=1-504\sum_{n=1}^\inf \bigl(\sum_{d|n} d^5\bigr)q^n
= 1 - 504q - 16632q^2 - 122976q^3 -\cdots$$ be the Eisenstein
series of weight $4$ and $6$ on $\g$,
$$\D(\t)=\frac{E_4(\t)^3-E_6(\t)^2}{1728}= q - 24q^2 + 252q^3 -
1472q^4 +\cdots$$
the ``discriminant'' cusp form of weight $12$ and
$$j(\t)=\frac{E_4(\t)^3}{\D(\t)} = \frac1q + 744 + 196884q +
21493760q^2 + \cdots$$
the modular  invariant.  For an integer $N$,
let $\G_0(N)$ denote the modular group of level $N$ defined by
$$\G_0(N)=\left\{\left(\matrix a&b\\c&d \endmatrix\right)\in\g  |\;
  c\equiv 0\mod N\right\}.$$
Define
$$\align  E_2^{(2)}(\t) &:= 2E_2(2\t)-E_2(\t)\\&=1+24\sum_{n=1}^\inf
\bigl(
\sum_{d|n\atop d:\text{odd}} d\; \bigr)q^n= 1+24q+24q^2+96q^3+\cdots,\\
\D_4^{(2)}(\t) &:= \frac{\eta(2\t)^{16}}{\eta(\t)^8}= \sum_{n=1}^\inf
\bigl(\sum_{d|n\atop d:\text{odd}}
(n/d)^3\bigr) q^n=q+8q^2+28q^3+64q^4+\cdots,\\
j^{(2)}(\t)
&:=\frac{E_2^{(2)}(\t)^2}{\D_4^{(2)}(\t)}=\frac1q+40+276q-2048q^2+\cdots,
\endalign $$
where $$\eta(\t)=q^{\frac1{24}}\prod_{n=1}^\inf (1-q^n)
=q^\frac{1}{24} - q^\frac{25}{24}- q^\frac{49}{24} + q^\frac{121}{24}+
\cdots$$
is the Dedekind eta function.  The functions $E_2^{(2)}(\t)$
and $\D_4^{(2)}(\t)$ are modular forms of respective weights $2$ and
$4$ on the group $\G_0(2)$, and $j^{(2)}(\t)$ is a $\G_0(2)$-invariant
function which generates the field of modular functions on $\G_0(2)$
(the normalized function $j^{(2)}(\t)-40$ is often referred to as the
``Hauptmodul'' for the group $\G_0(2)$).  In addition, the function
$$\sqrt{\D_4^{(2)}(\t)}=\frac{\eta(2\t)^{8}}{\eta(\t)^4}=
\sum_{n=1\atop n:\text{odd}}^\inf \bigl(\sum_{d|n} d \bigr)
q^{\frac{n}2}= q^{\frac12}+4q^{\frac32}+6q^\frac52+8q^\frac72+\cdots$$
will also appear in the formula. This is of weight $2$ belonging to
the principal congruence subgroup $$\Gamma(2)=\left\{\left(\matrix
    a&b\\c&d \endmatrix\right)\in\g  |\; a\equiv d\equiv 1,\ b\equiv
  c\equiv 0\mod 2\right\},$$
which is of index $2$ in
$\Gamma_0(2)$. (This is seen from the transformation formula of the
eta function.)

Similarly, we define
$$\align  E_1^{(3)}(\t) &:= 1+6\sum_{n=1}^\inf
\bigl(\sum_{d|n} \left(\frac{d}{3}\right)\bigr) q^n= 1+6q+6q^3+6q^4+\cdots,\\
\D_3^{(3)}(\t) &:= \frac{\eta(3\t)^{9}}{\eta(\t)^3}=\sum_{n=1}^\inf
\bigl(\sum_{d|n} \left(\frac{d}{3}\right)(n/d)^2\bigr) q^n
=q+3q^2+9q^3+13q^4+\cdots,\\
j^{(3)}(\t)
&:=\frac{E_1^{(3)}(\t)^3}{\D_3^{(3)}(\t)}=\frac1q+15+54q-76q^2-\cdots,
\endalign $$
where $\left(\dfrac{d}{3}\right)$ is the Legendre
character.  Both $E_1^{(3)}$ and $\D_3^{(3)}(\t)$ are modular forms of
weights $1$ and $3$, respectively, with character
$\left(\dfrac{d}{3}\right)$ on the group $\Gamma_0(3)$, while the
function $j^{(3)}(\t)$ is a generator of the field of modular
functions on $\Gamma_0(3)$.  We also need
$$\bigl(\D_3^{(3)}(\t)\bigr)^{\frac13} =
\frac{\eta(3\t)^{3}}{\eta(\t)}=\sum_{n=1\atop n\not\equiv 0 (3)}^\inf
\bigl(\sum_{d|n} \left(\frac{d}{3}\right)\bigr) q^n
=q^\frac13+q^\frac43+2q^\frac73+2q^\frac{13}3+\cdots,$$
which, as
follows from the transformation formula of the eta function (as in the
case of $\sqrt{\D_4^{(2)}(\t)}$), is of weight $1$ (with the same
character) on the subgroup
$$\Gamma_0^0(3)=\left\{\left(\matrix a&b\\c&d
    \endmatrix\right)\in\g  |\; b\equiv c\equiv 0\mod 3\right\}$$
of index $3$ of $\Gamma_0(3)$. Finally, let
$$\align  E_2^{(4)}(\t) &:= \frac13\left(4E_2(4\t)-E_2(\t)\right)
=1 + 8q + 24q^2 + 32q^3 +24q^4 +\cdots,\\
\D_2^{(4)}(\t) &:= \frac{\eta(4\t)^{8}}{\eta(2\t)^4} =\sum_{n=1\atop
  n:\text{odd}}^\inf \bigl(\sum_{d|n}d\bigr)q^n
=q+4q^3+6q^5+8q^7+\cdots,\\ j^{(4)}(\t)
&:=\frac{E_2^{(4)}(\t)}{\D_2^{(4)}(\t)}=\frac1q+8+20q-62q^3+216q^5-\cdots.
\endalign $$
The functions $E_2^{(4)}(\t)$ and $\D_2^{(4)}(\t)$ are
both modular of weight $2$ on $\Gamma_0(4)$ and $j^{(4)}(\t)$
generates the field of modular functions on $\Gamma_0(4)$. Further we need
$$\bigl(E_2^{(4)}(\t)\bigr)^{\frac14}=\theta_3(2\t)=\sum_{n\in\Z}q^{n^2}
=1+2q+2q^4+2q^9+\cdots$$
which is of weight $1/2$ on $\Gamma_0(4)$ and
$$\bigl(\D_2^{(4)}(\t)\bigr)^{\frac14}=\frac{\eta(4\t)^{2}}{\eta(2\t)}
=\frac12\theta_2(2\t)=\frac12\sum_{n\in\Z}q^{(n+\frac12)^2}=q^\frac14
+q^\frac94+q^\frac{25}4+q^\frac{49}4+\cdots$$
which is of weight $1/2$
on $$\Gamma_0^0(4)=\left\{\left(\matrix a&b\\c&d
    \endmatrix\right)\in\g  |\; b\equiv c\equiv 0\mod 4\right\}.$$
Here, $\theta_3$ and $\theta_2$ are Jacobi's theta functions (we do
not use this notation later on). That the $\theta_2(2\tau)$ belongs
to $\Gamma_0^0(4)$ is seen from the transformation formula for theta
functions or, alternatively, from the fact that
$\theta_2(8\t)=\theta_3(2\t)-\theta_3(8\t)$ is a modular form of
weight $1/2$ on $\Gamma_0(16)$ and that
$$\Gamma_0^0(4)=\left(\matrix 4&0\\0&1
    \endmatrix\right)\Gamma_0(16)\left(\matrix \frac14 &0\\0&1
    \endmatrix\right).$$

Recall the definition of the Gauss hypergeometric
series $F={}_2F_1$: $$F(a,b,c;x)=\sum_{n=0}^\inf \frac{(a)_n(b)_n}
{(c)_n}\frac{x^n}{n!},$$
where $(a)_n$ denotes $a(a+1)\cdots (a+n-1)$.

\proclaim{Theorem 1} \roster \item"{(i)}" When $k$ is an even integer
congruent to $0$ or $4$ modulo $6$, the equation \eq  has a one
dimensional space of solutions which are holomorphic modular forms of
weight $k$ on $\g$, a generator of which is given by
$$\align & E_4(\tau)^{\frac{k}4}
F(-\frac{k}{12},-\frac{k-4}{12},-\frac{k-5}6;\frac{1728}{j(\t)}) \\
&=\sum_{0\le i\le k/12}\frac{(-\frac{k}{12})_i(-\frac{k-4}{12})_i}
{(-\frac{k-5}6)_i\,i!}1728^i\D(\t)^iE_4(\t)^{\frac{k}4-3i}=1+\O(q)
\endalign $$
if $k\equiv 0,4 \mod 12$, and by
$$\align &E_4(\tau)^{\frac{k-6}4}E_6(\tau)
F(-\frac{k-6}{12},-\frac{k-10}{12},-\frac{k-5}6;\frac{1728}{j(\t)})\\
&=E_6(\tau)\sum_{0\le i\le
  (k-6)/12}\frac{(-\frac{k-6}{12})_i(-\frac{k-10}{12})_i}
{(-\frac{k-5}6)_i \,i!}1728^i\D(\t)^iE_4(\t)^{\frac{k-6}4-3i}\\
&=1+\O(q) \endalign $$
if $k\equiv 6,10 \mod 12$.
\item"{(ii)}" When $k$ is an even integer congruent to $2$ modulo $6$, \eq has
  a two dimensional space of modular solutions$:$ The function
  $$\align
  &E_2^{(2)}(\t)^{\frac{k}2}F(-\frac{k}4,-\frac{k-2}4,-\frac{k-5}6;
  \frac{64}{j^{(2)}(\t)})\\
  &=\sum_{0\le i\le k/4}\frac{(-\frac{k}{4})_i(-\frac{k-2}{4})_i}
  {(-\frac{k-5}6)_i
    \,i!}64^i\D_4^{(2)}(\t)^iE_2^{(2)}(\t)^{\frac{k}2-2i}
  =1+\O(q)\endalign $$
  which is a holomorphic modular form of  weight
  $k$ on $\G_0(2)$ generates a one dimensional subspace and the
  function
  $$\align &\D_4^{(2)}(\t)^{\frac{k+1}6}E_2^{(2)}(\t)^{\frac{k-2}6}
  F(-\frac{k-2}{12},-\frac{k-8}{12},\frac{k+7}6;\frac{64}{j^{(2)}(\t)})\\
  &=\sum_{0\le i\le
    (k-2)/12}\frac{(-\frac{k-2}{12})_i(-\frac{k-8}{12})_i}
  {(\frac{k+7}6)_i \,i!}64^i\D_4^{(2)}(\t)^{\frac{k+1}6+i}
  E_2^{(2)}(\t)^{\frac{k-2}6-2i}\\
  &=q^{\frac{k+1}6}+\O(q^{\frac{k+7}6})\endalign$$
  which is of weight
  $k$ on $\Gamma(2)$ constitutes the other generator.
\item"{(iii)}" When $k$ is an odd integer congruent to $1$ or $3$ modulo $6$,
  \eq has a two dimensional space of modular solutions$:$ The function 
  $$\align &E_1^{(3)}(\t)^kF(-\frac{k}3,-\frac{k-1}3,-\frac{k-5}6;
  \frac{27}{j^{(3)}(\t)})\\
  &=\sum_{0\le i\le k/3}\frac{(-\frac{k}{3})_i(-\frac{k-1}{3})_i}
  {(-\frac{k-5}6)_i \,i!}27^i\D_3^{(3)}(\t)^iE_1^{(3)}(\t)^{k-3i}
  =1+\text{O}(q)\endalign $$
  which is a holomorphic modular form of
  weight $k$ on $\G_0(3)$ generates a one dimensional subspace and the
  function $$\align
  &\D_3^{(3)}(\t)^{\frac{k+1}6}E_1^{(3)}(\t)^{\frac{k-1}2}
  F(-\frac{k-1}{6},-\frac{k-3}{6},\frac{k+7}6;\frac{27}{j^{(3)}(\t)})\\
  &=\sum_{0\le i\le
    (k-1)/6}\frac{(-\frac{k-1}{6})_i(-\frac{k-3}{6})_i}
  {(\frac{k+7}6)_i \,i!}27^i\D_3^{(3)}(\t)^{\frac{k+1}6+i}
  E_1^{(3)}(\t)^{\frac{k-1}2-3i}\\
  &=q^{\frac{k+1}6}+\O(q^{\frac{k+7}6})\endalign$$
  which is of weight
  $k$ on $\Gamma_0^0(3)$ constitutes the other generator.
\item"{(iv)}" When $k$ is half an integer congruent to $\frac12$ modulo $3$,
  \eq has a two dimensional space of modular solutions$:$ The function
  $$\align
  &E_2^{(4)}(\t)^{\frac{k}2}F(-\frac{2k-1}6,-\frac{k}2,-\frac{k-5}6;
  \frac{16}{j^{(4)}(\t)})\\
  &=\sum_{0\le i\le
    (2k-1)/6}\frac{(-\frac{2k-1}{6})_i(-\frac{k}{2})_i}
  {(-\frac{k-5}6)_i
    \,i!}16^i\D_2^{(4)}(\t)^iE_2^{(4)}(\t)^{\frac{k}2-i}
  =1+\text{O}(q)\endalign $$
  which is a holomorphic modular form of
  half-integral weight $k$ on $\G_0(4)$ generates a one dimensional
  subspace and the function $$\align
  &\D_2^{(4)}(\t)^{\frac{k+1}6}E_2^{(4)}(\t)^{\frac{2k-1}6}
  F(-\frac{2k-1}{6},-\frac{k-2}{6},\frac{k+7}6;\frac{16}{j^{(4)}(\t)})\\
  &=\sum_{0\le i\le
    (2k-1)/6}\frac{(-\frac{2k-1}{6})_i(-\frac{k-2}{6})_i}
  {(\frac{k+7}6)_i \,i!}16^i\D_2^{(4)}(\t)^{\frac{k+1}6+i}
  E_2^{(4)}(\t)^{\frac{2k-1}6-i}\\
  &=q^{\frac{k+1}6}+\O(q^{\frac{k+7}6})\endalign$$
  which is of weight
  $k$ on $\Gamma_0^0(4)$ constitutes the other generator.  \endroster
  \endproclaim

\example{Remark} 1. The case (i) is contained in [1, Theorem 5].

2. We expect by numerical evidence no other modular (on a congruence
subgroup) solution to $(\#)_k$, at least when $k$ is an integer or half an
integer (denominators of coefficients of power series solution in
other cases seem not to be bounded). We however have no proof of this
speculation.    \endexample

Before giving a proof, we introduce the operator $\p_k$ defined by
$$\p_k(f)(\t)=f'(\t)-\frac{k}{12}E_2(\t)f(\t).$$
By the quasimodular property of $E_2(\t)$ which reads
$$E_2\left(\frac{a\t+b}{c\t+d}\right)=(c\t+d)^2E_2(\t)+\frac6{\pi i}
c(c\t+d) \qquad (\pmatrix a & b \\ c & d \endpmatrix \in\g),$$
we see
that if $f$ is modular of weight $k$ on a subgroup of $\g$,  then
$\p_k(f)$  is modular of weight $k+2$ on the same group.  If $f$ and
$g$ have weights $k$ and $l$, the Leibniz rule
$$\p_{k+l}(fg)=\p_k(f)g+f\p_l(g)$$
holds. We shall often drop the
suffix  of the operator $\p_k$ since the weights of modular forms we
shall be considering are clear in most cases. With this operator, the
equation \eq can be written as
$$(\#')_k\qquad\p_{k+2}\p_k(f)(\t)=\frac{k(k+2)}{144}E_4(\t)f(\t),
\hskip60pt$$ (use $E_2'=(E_2^2-E_4)/12$).

\demo{Proof of Theorem 1} Given a specific modular form in terms of
known forms, it is a straightforward task to check if it satisfies \eq
or not.  We give a proof of the first part of (ii) to illustrate the calculation, the
remaining cases being similar.  Write $A=\D_4^{(2)}, B=E_2^{(2)}$ to
ease notation.  We have $$\p_4(A)=\frac23AB,\
\p_2(B)=32A-\frac16B^2,$$
and $$E_4=192A+B^2.$$
(To establish these
kinds of identities, it is enough to check that the first several
Fourier coefficients coincide, since both sides of these equations are
holomorphic modular forms of weight $6,4$ and $4$ on $\G_0(2)$.)  Using
the first two we obtain
$$\p^2(A^iB^{\frac{k}2-2i})
=aA^iB^{\frac{k}2-2i+2}+bA^{i+1}B^{\frac{k}2-2i}
+cA^{i+2}B^{\frac{k}2-2i-2}$$
with $$a=\frac{(k-12i)(k-12i+2)}{144},
b=-\frac83((k-4i)(k-12i-4)-8i), c=256(k-4i)(k-4i-2).$$Hence, for
$$f=\sum_{0\le i\le k/4}c_i A^iB^{\frac{k}2-2i}\quad\text{with} \quad
c_i=64^i\frac{(-\frac{k}4)_i(-\frac{k-2}4)_i}{(-\frac{k-5}6)_i \,i!},$$
we  have
$$\p^2(f)-\frac{k(k+2)}{144}E_4f =\sum_{0\le i\le k/4}c_i'
A^{i+1}B^{\frac{k}2-2i},$$
where
$$\align c_i'&=\frac{(k-12i-12)(k-12i-10)}{144}c_{i+1}
-\frac83((k-4i)(k-12i-4)-8i)c_i\\
&+256(k-4i+4)(k-4i+2)c_{i-1}-\frac{k(k+2)}{144}(c_{i+1}+192c_i)
\\
&=c_i\times
\biggl\{-\frac16\frac{(k-12i-12)(k-12i-10)(k-4i)(k-4i-2)}
{(k-6i-5)(i+1)}\\
  &-\frac83((k-4i)(k-12i-4)-8i)-\frac{64}{6} (k-6i+1)i\\
  &-\frac{k(k+2)}{144}\left(-24\frac{(k-4i)(k-4i-2)}{(k-6i-5)(i+1)}
    +192\right)\biggr\}, 
\endalign$$
which turns out to be identically
$0$.\endpf

The next theorem describes a solution in the case of $k\equiv 5 \mod 6$.
Here we come across a different phenomenon: the equation \eq has {\it
  quasimodular} solutions of weight $k+1$ rather than $k$. Recall that
an element of degree $k$ in the graded ring
$\C[E_2(\t),E_4(\t),E_6(\t)]$, where the generators $E_2, E_4, E_6$
have degrees $2, 4$, and $6$, is referred to as a  quasimodular form
of weight $k$ (on $\g$).  Define a sequence of polynomials $P_n(x)\
(n=0,1,2,\dots)$ by
$$P_0(x)=1,\  P_1(x)=x,\  \ P_{n+1}(x)=xP_n(x)+\l_nP_{n-1}(x)\
(n=1,2,\dots)$$
where $$\l_n=12\frac{(6n+1)(6n+5)}{n(n+1)}.$$
First
few examples are $$P_2(x)=x^2+462, P_3(x)=x^3+904x,
P_4(x)=x^4+1341x^2+201894,\dots.$$
Clearly $P_n(x)$ is an even or odd
polynomial according as $n$ is even or odd. We also define a series of
``companion'' polynomials $Q_n(x)$ by the same recursion  (with
different initial values):
$$Q_0(x)=0,\  Q_1(x)=1,\ \ Q_{n+1}(x)=xQ_n(x)+\l_nQ_{n-1}(x)\
(n=1,2,\dots),$$
a few examples being $$Q_2(x)=x, Q_3(x)=x^2+442,
Q_4(x)= x^3+879x,\dots.$$
$Q_n(x)$ has the opposite parity: It is even
if $n$  is odd and odd if $n$ is even.

\proclaim{Theorem 2} Let $k=6n+5\ (n=0,1,2,\dots)$. The following 
quasimodular form of weight $k+1$ on $\g$ is a solution of \eq$:$
$$\sqrt{\D(\t)}^nP_n\bigl(\frac{E_6(\t)}{\sqrt{\D(\t)}}\bigr)\frac{E_4'(\t)}{240}
-\sqrt{\D(\t)}^{n+1}Q_n\bigl(\frac{E_6(\t)}{\sqrt{\D(\t)}}\bigr).$$
\endproclaim

\example{Remark} Because of the parities of $P_n(x)$ and $Q_n(x)$, the
appearance of $\sqrt{\D(\t)}$ in the formula is in fact superficial,
i.e., it cancels out and the formula gives an element in $\Q[
E_2(\t),E_4(\t),E_6(\t)]$ (note $E_4'(\t)=(E_2(\t)E_4(\t)-E_6(\t))/3$). 
\endexample
 
Since the polynomials $P_n(x)$ and $Q_n(x)$ are {\it not}
hypergeometric polynomials, it is not at all straightforward to show
the expression in the  theorem satisfies $(\#)_k$. We shall give a proof of
the theorem in the next section where an inductive  structure of
solutions of \eq with varying $k$ is discussed.  \medskip

A few more words about the polynomials $P_n(x)$ and $Q_n(x)$: If we
replace the constant $\l_n=12(6+1/n)(6-1/(n+1))$ in the recursion of
$P_n$ and $Q_n$ by $12(6+(-1)^n/(n-1))(6+(-1)^n/n)$, (the even part
of) the resulting polynomials give ``Atkin's orthogonal polynomials''
(see [1, Sections 4 and 5]). In our case, the polynomials $P_n(x)$ and
$Q_n(x)$ are connected to the convergents of the continued fraction
expansion of
$$\frac1{240}\frac{E_4'}{E_6}=
\frac1{j}\frac{F(\frac{13}{12},\frac{17}{12},2;\frac{1728}j)}
{F(\frac{1}{12},\frac{5}{12},1;\frac{1728}j)}
=\frac1j+\frac{1266}{j^2}+\frac{1806960}{j^3}+\cdots, $$
in a manner similar to Atkin's case where the function involved was
$$\frac{E_2E_4}{E_6j}=
\frac1{j}\frac{F(\frac{13}{12},\frac{5}{12},1;\frac{1728}j)}
{F(\frac{1}{12},\frac{5}{12},1;\frac{1728}j)}
=\frac1j+\frac{720}{j^2}+\frac{911520}{j^3}+\cdots. $$
It should be
possible, by a method analogous to that in [1], to establish
properties like a closed formula or a differential equation for our
$P_n$ and $Q_n$.  We however do not pursue this here.

\section{\S3. Inductive structure of solutions}

For modular forms $f(\t)$ and $g(\t)$ of weights $k$ and  $l$,
define a modular form $[f(\t), g(\t)]$ of weight $k+l+2$ by
$$[f(\t), g(\t)]=kf(\t)g'(\t)-lf'(\t)g(\t)$$ 
(``Rankin-Cohen bracket'' of degree $1$). The right-hand 
side may also be written as $$kf(\t)\p_l(g(\t))-l\p_k(f(\t))g(\t).$$ 
 
\proclaim{Lemma} Suppose $f=F_k(\t)$ satisfies  $(\#')_k$. Then we have
$$\p_{k+6}([F_k(\t),E_4(\t)])=\frac{k-4}{18}[F_k(\t),E_6(\t)]$$ and
$$\p_{k+8}([F_k(\t),E_6(\t)])=\frac{k-6}{8}E_4(\t)[F_k(\t),E_4(\t)].$$
Here, in the Rankin-Cohen brackets, the function $F_k(\t)$ is regarded
as being of weight $k$.
\endproclaim

\pf Using $\p(E_4)=-E_6/3$ and $\p(E_6)=-E_4^2/2$, we
have
$$\align \p([F_k, E_4])&=\p\bigl(-\frac{k}3F_kE_6-4\p(F_k)E_4\bigr)\\
&=-\frac{k}3\p(F_k)E_6+\frac{k}6F_kE_4^2-4\biggl(\frac{k(k+2)}{144}E_4F_k\biggr)
E_4+\frac43\p(F_k)E_6\\
&=-\frac{k(k-4)}{36}F_kE_4^2-\frac{k-4}3\p(F_k)E_6\\
&=\frac{k-4}{18}[F_k,E_6]. \endalign $$ The other identity is shown
similarly. \endpf

\proclaim {Proposition 1} \roster \item"{(i)}" If $F_k(\t)$ is a solution of
$(\#)_k$, then $[F_k(\t), E_4(\t)]/\D(\t)$ is a solution of \eqm.
\item"{(ii)}" Assume $k\ne 0,4,5$. Given a solution $F_k(\t)$ of $(\#)_k$, put 
$$F_{k-6}(\t)=\frac{k-5}{288k(k-4)}\dfrac{[F_k(\t), E_4(\t)]}{\D(\t)}$$
and define $F_{k+6i}(\t)\ (i=1,2,3,\dots)$ successively by the recursion
$$F_{k+6i+6}(\t)=E_6(\t)F_{k+6i}(\t)+\mu_i^{(k)}\D(\t)F_{k+6i-6}(\t) \quad
(i=0,1,2,\dots)$$
where
$$\mu_i^{(k)}=432\frac{(k+6i)(k+6i-4)}{(k+6i+1)(k+6i-5)}.$$
Then
$F_{k+6i}(\t)$  is a solution of $(\#)_{k+6i}$ for every $i$.
\endroster \endproclaim

\pf Noting $\p(\D(\t))=0$ (since $\D'=E_2\D$), the assertion (i) is
readily shown as  a direct consequence of Lemma. In fact, we have
$$\p^2\biggl(\frac{[F_k,E_4]}{\D}\biggr)=\frac{\p^2([F_k,E_4])}{\D}
=\frac{k-4}{18}\frac{\p([F_k,E_6])}{\D}
=\frac{(k-6)(k-4)}{144}E_4\frac{[F_k,E_4]}{\D},$$
which shows  that the function $[F_k,
E_4]/\D$ satisfies the equivalent form $(\#')_{k-6}$ of $(\#)_{k-6}$.

For (ii), we first show the following.

\proclaim{Lemma} If $G_k$ and $G_{k-6}$ are solutions of  \eq and \eqm
respectively, then the function $G_{k+6}:=E_6G_k+\D G_{k-6}$  is a
solution of \eqp if and only if the relation
 $$[G_k, E_4]=\frac23(k+1)\D G_{k-6}$$ holds. \endproclaim
 
 \pf Using $$\align \p(E_6G_k)&=-\frac{E_4^2}2G_k+E_6\p(G_k)
 =-\frac16[G_k,E_6]-\frac{k+6}{12}E_4^2G_k,\\
 \p(E_4^2G_k)&=-\frac23E_4E_6G_k+E_4^2\p(G_k)
 =-\frac14[G_k,E_4]-\frac{k+4}{12}E_6G_k,\endalign$$
 and the preceding
 lemma, we have
 $$\align \p^2(G_{k+6})&=\p^2(E_6G_k)+\D\p^2(G_{k-6})\\
 &=\p\biggl(-\frac16[G_k,E_6]-\frac{k+6}{12}E_4^2G_k\biggr)+\D\p^2(G_{k-6})\\
 &=
 -\frac16\frac{k-6}8E_4[G_k,E_4]-\frac{k+6}{12}\bigl(-\frac14E_4[G_k,E_4]
 -\frac{k+8}{12}E_4E_6G_k\bigr)\\
 &+\D\frac{(k-6)(k-4)}{144}E_4G_{f-6}\\
 &=\frac14E_4[G_k,E_4]+\frac{(k+6)(k+8)}{144}E_4E_6G_k+
 \frac{(k-6)(k-4)}{144}E_4\D G_{k-6}.\endalign $$
 Hence we obtain
$$\p^2(G_{k+6})-\frac{(k+6)(k+8)}{144}E_4G_{k+6}=\frac{E_4}4\left([G_k,E_4]
-\frac23(k+1)\D G_{k-6}\right).$$ The lemma follows from this. \endpf 

We return to the proof of (ii) of Proposition 1. By (i), the  function
$\mu_0F_{k-6}$ is a solution of \eqm and it is so defined that the
relation
$$
[F_k,E_4]=\frac23(k+1)\D(\mu_0 F_{k-6})$$
holds. By the lemma we
have just shown, we conclude that $F_{k+6}=E_6F_k+\mu_0\D F_{k-6}$ is
a solution of \eqp. Moreover, the pair $F_{k+6}$ and $\mu_1F_k$ satisfy
$$[F_{k+6},E_4]=\frac23(k+7)\D (\mu_1F_k)$$ as the following calculation
shows;
$$\align [F_{k+6},E_4]&=[E_6F_k+\mu_0\D F_{k-6},E_4]
=[E_6F_k+\frac3{2(k+1)}[F_{k},E_4],E_4]\\
&=(k+6)E_6F_k\p(E_4)-4\p(E_6F_k)E_4\\
&+\frac3{2(k+1)}\left((k+6)[F_k,E_4]\p(E_4)-4\p([F_k,E_4])E_4\right)\\
&=\frac{k+6}3E_6^2F_k+2E_4^3F_k-4E_4E_6\p(F_k)\\
&-\frac{k+6}{2(k+1)}E_6\left(-\frac{k}3E_6F_k-4E_4\p(F_k)\right)\\
&-\frac6{k+1}E_4\frac{k-4}{18}\left(-\frac{k}2 E_4^2F_k-6E_6\p(F_k)\right)\\
&=\frac{(k+2)(k+6)}{6(k+1)}(E_4^3-E_6^2)F_k =\frac23(k+7)\mu_1\D
F_k. \endalign$$
Thus, applying the lemma again, we conclude $F_{k+12}$
is a solution of $(\#)_{k+12}$. By replacing $k$ with $k+6$ in the
above  calculation, we see this procedure continues inductively and
the proposition is proved. \endpf

\demo{Proof of Theorem 2} For $k=6n+5\ (n=0,1,2,\dots)$, denote by $F_k(\t)$
the form in the theorem,
$$F_k(\t)=\sqrt{\D(\t)}^nP_n\bigl(\frac{E_6(\t)}{\sqrt{\D(\t)}}\bigr)
\frac{E_4'(\t)}{240}
-\sqrt{\D(\t)}^{n+1}Q_n\bigl(\frac{E_6(\t)}{\sqrt{\D(\t)}}\bigr).$$
We see directly that $F_{11}=E_6E_4'/240-\D$ is a solution of
$(\#)_{11}$ and that
$$F_5=\frac{E_4'}{240}=\frac{11-6}{288\cdot
  11\cdot(11-4)}\frac{[F_{11},E_4]}{\D}.$$
By Proposition 1 (ii), it is
then enough to show that the $F_k$'s satisfy the recursion in
Proposition 1 (ii). Alternatively, we may start with checking that
$F_5=E_4'/240$ satisfies $(\#)_5$ and calculate $[F_5,E_4]=-4\D$. By
this, if we put $F_{-1}=1$ and $\mu_0^{(5)}=-1$, then the initial
condition of Proposition 1 (ii) is satisfied and also we have
$F_{11}=E_6F_5+\mu_0^{(5)}\D F_{-1}$, so again the task is that the
$F_k$'s satisfy the recursion in Proposition 1 (ii). Noticing the
constant $\mu_i^{(k)}$ in Proposition 1 is equal to
$$432\frac{(6n+5+6i)(6n+5+6i-4)}{(6n+5+6i+1)(6n+5+6i-5)}=
12\frac{(6n+6i+5)(6n+6i+1)}{(n+i+1)(n+i)}=\l_{n+i},$$ we
have, by the recursion of $P_n$ and $Q_n$,
$$
\align &E_6F_{k+6i}+\mu_i^{(k)}\D F_{k+6i-6}\\
&=E_6\biggl(\sqrt{\D}^{n+i}P_{n+i}\bigl(\frac{E_6}{\sqrt{\D}}\bigr)\frac{E_4'}{240}
-\sqrt{\D}^{n+i+1}Q_{n+i}\bigl(\frac{E_6}{\sqrt{\D}}\bigr)\biggr)\\
&+\l_{n+i}\D\biggl(
\sqrt{\D}^{n+i-1}P_{n+i-1}\bigl(\frac{E_6}{\sqrt{\D}}\bigr)\frac{E_4'}{240}
-\sqrt{\D}^{n+i}Q_{n+i-1}\bigl(\frac{E_6}{\sqrt{\D}}\bigr)\biggr)\\
&=\sqrt{\D}^{n+i+1}\biggl(\frac{E_6}{\sqrt{\D}}
P_{n+i}\bigl(\frac{E_6}{\sqrt{\D}}\bigr)+\l_{n+i}
P_{n+i-1}\bigl(\frac{E_6}{\sqrt{\D}}\bigr)\biggr)\frac{E_4'}{240}\\
&-\sqrt{\D}^{n+i+2}\biggl(\frac{E_6}{\sqrt{\D}}
Q_{n+i}\bigl(\frac{E_6}{\sqrt{\D}}\bigr)+\l_{n+i}
Q_{n+i-1}\bigl(\frac{E_6}{\sqrt{\D}}\bigr)\biggr)\\
&=\sqrt{\D}^{n+i+1}P_{n+i+1}\bigl(\frac{E_6}{\sqrt{\D}}\bigr)\frac{E_4'}{240}
-\sqrt{\D}^{n+i+2}Q_{n+i+1}\bigl(\frac{E_6}{\sqrt{\D}}\bigr)\\
&=F_{k+6i+6}. \endalign $$
We therefore conclude that $F_{k+6i}$ is a
solution of $(\#)_{k+6i}$ for every $i$. \endpf

\example{Remark} The recursion in Proposition 1 is also satisfied by
the  modular solutions given in Theorem 1. In that case, the recursion
reduces to relations of hypergeometric series but we should point out
that those relations are not in general the relations referred to as
the contiguous relations of Gauss.  For instance, the first solution
in (ii) of Theorem 1 (the case of $k\equiv 2 \bmod 6$) gives an identity
of hypergeometric series of the form (which is true for any $k$)
$$\align & x^3 F\bigl(-\frac{k+6}4,-\frac{k+4}4,-\frac{k+1}6;\frac{64}x\bigr)\\
&=(x^3-576x^2) F\bigl(-\frac{k}4,-\frac{k-2}4,-\frac{k-5}6;\frac{64}x\bigr)\\
&+432\frac{k(k-4)}{(k-5)(k+1)}(x^2-128x+4096)
F\bigl(-\frac{k-6}4,-\frac{k-8}4,-\frac{k-11}6;\frac{64}x\bigr),
\endalign $$
where $x=j^{(2)}=B^2/A$ (the polynomials on the
right-hand side come from the identities
$$E_6=B^3-576AB,\quad \D=AB^4-128A^2B^2+4096 A^3,$$
notation being as
in the proof of Theorem 1).  Instead of checking  these sorts of
hypergeometric identities case by case in order to  show the recursion
for solutions in Theorem 1, we argue as follows.  By the form of the
recursion in Proposition 1 (ii), if $F_k$ has $q$-expansion of the
form $1+\O(q)$ (with only integral powers of $q$), so does every
$F_{k+6i}$ (use $E_6=1+\O(q)$ and $\D=q+\O(q^2)$).  Hence, starting
with a solution in Theorem 1 of the form $1+\O(q)$ with $0<k\le 6$ or
$k=10$, we conclude by the uniqueness of solution of  this form that
the recursively determined solutions must coincide the ones in Theorem
1. The case of other solutions with $q^{(k+1)/6}+\O(q^{(k+7)/6})$ is
similar.  Hence, the above hypergeometric identity or the other
corresponding  ones may be regarded as consequences of Theorem
1. \endexample

Proposition 1 shows that if we find any solution of $(\#)_k$, we can
construct  solutions of \eq for larger $k$ in the same residue class
modulo $6$.  Conversely, any solution ``comes from'' lower ones in
this way. In fact, suppose $F_k$ is a solution of $(\#)_k$. Applying
Proposition 1 (i) twice we see not only the function $[F_k,E_4]/\D$
is a solution of \eqm but also the function $[[F_k,E_4]/\D,E_4]/\D$ is
a solution of $(\#)_{k-12}$. Put
$$F_{k-6}=\frac{k-5}{288k(k-4)}\frac{[F_k, E_4]}{\D}$$
and
$$F_{k-12}=\frac{k-11}{288(k-6)(k-10)}\frac{[F_{k-6}, E_4]}{\D}.$$
With these we have
$$F_k=E_6F_{k-6}+\mu_0^{(k-6)}F_{k-12}.$$
This tells us that the $F_k$
is obtained from $F_{k-6}$ as in Proposition 1 (ii).  Taking account of
this consideration and the fact that if $F_k$ is a  holomorphic
modular form then so is $[F_k,E_4]/\D$, in order to find a modular
solution of $(\#)_k$, we may restrict ourselves to find one in the range
$0\le k <6$. (Since $[1,E_4]=[E_4,E_4]=0$, Proposition 1 (ii) is valid
for $k=0$ or $4$ if we take $F_0=1$ or $F_4=E_4$ with $F_{-6}=0$ or
$F_{-2}=0$.)

Remarkably enough, in a completely different context (deformation of
singularities etc.) the equation \eq (in its yet another equivalent
form different from $(\#')_k$) for  integral or half-integral $k$ in this
range was considered by Ikuo Satake (private communication, related
works are [2,3]).  There, the solutions are expressed as theta
functions associated with positive definite lattices, a table of which
is given below:
$$
\vbox{\offinterlineskip \halign{\strut \hfil$#$\hfil\enskip  &
    \vrule width 0.8pt $#$& \enskip \hfil$#$\enskip &\vrule$#$&\enskip
    \hfil$#$\enskip &\vrule$#$&\enskip\hfil$#$\enskip
    &\vrule$#$&\enskip\hfil$#$\enskip
    &\vrule$#$&\enskip\hfil$#$\enskip
    &\vrule$#$&\enskip\hfil$#$\enskip
    &\vrule$#$&\enskip\hfil$#$\enskip
    &\vrule$#$&\enskip\hfil$#$\enskip
    &\vrule$#$&\enskip\hfil$#$\enskip
    &\vrule$#$&\enskip\hfil$#$\enskip
    &\vrule$#$&\enskip\hfil$#$\enskip
    &\vrule$#$&\enskip\hfil$#$\enskip \cr k && 0 && \frac12\enskip  &&
    1\,\enskip && \frac32\, && 2\,\enskip &&  \frac52\, && 3\,\enskip
    && \frac72\enskip && 4\,\enskip && \frac92\,  &&  5\; &&
    \frac{11}2\cr \noalign{\hrule height 0.8pt} F_k && 1 &&
    \Theta_{A_1} && \Theta_{A_2} && $---$ && \Theta_{D_4}  && $---$ &&
    \Theta_{E_6} && \Theta_{E_7} && \Theta_{E_8} && $---$ && $---$ &&
    $---$ \cr &&  && \Theta_{A_1^*} && \Theta_{A_2^*} &&  &&
    \Theta_{D_4^*} && &&  \Theta_{E_6^*} && \Theta_{E_7^*} && && &&
    &&\cr }} $$
Here, $\Theta_L$ is the theta function of a lattice
$L$ and  $A_n,D_n,E_n$ denote root lattices while $A_n^*,D_n^*,E_n^*$
are their duals.  In Satake's investigation, the existence of modular
solution in that range (e.g., we have a modular solution for $k=1/2$
but not for $k=3/2$) corresponds exactly to the existence of certain
affine root system. It seems to be an interesting question if the
existence or non-existence of modular solution of \eq is explained by
any number theoretical reason. Also, the meaning of the quasimodular
solution $E_4'$ (for $k=5$) and the ones for higher weights should be
clarified.

\section{\S4. Positiveness of Fourier coefficients}

As remarked in the previous section, some solutions of \eq for small
$k$ are theta series of positive definite lattices and  hence have
positive Fourier coefficients. For general $k$, we prove the following.

\proclaim{Theorem 3} All the solutions given in Theorem 1 with
$q$-expansion of the  form $q^{(k+1)/6}+\O(q^{(k+7)/6})$ have positive
Fourier coefficients.  \endproclaim

\pf First we prove the case $k\equiv 2\mod 6$. The solution
in question is
$$\sum_{0\le i\le (k-2)/12}\frac{(-\frac{k-2}{12})_i(-\frac{k-8}{12})_i}
{(\frac{k+7}6)_i \,i!}64^i\D_4^{(2)}(\t)^{\frac{k+1}6+i}
E_2^{(2)}(\t)^{\frac{k-2}6-2i}.$$
Since the product
$(-(k-2)/12)_i(-(k-8)/12)_i$ in the numerator in the coefficient is always
positive for $0\le i\le (k-2)/12$ and the forms
$$\align E_2^{(2)}(\t)&=1+24\sum_{n=1}^\inf \bigl(\sum_{d|n\atop
  d:\text{odd}}
d\;\bigr) q^n,\\
\D_4^{(2)}(\t)&=\sum_{n=1}^\inf \bigl(\sum_{d|n\atop d:\text{odd}}
(n/d)^3\bigr) q^n \endalign $$
have positive coefficients, the
solution above clearly has positive Fourier coefficients.  The cases
$k\equiv 1,3\mod 6$ are treated similarly:  The product
$(-(k-1)/6)_i(-(k-3)/6)_i$ in the formula is positive and the forms
$$E_1^{(3)}(\t) = 1+6\sum_{n=1}^\inf
\bigl(\sum_{d|n} \left(\frac{d}{3}\right)\bigr) q^n$$ and
$$\D_3^{(3)}(\t) = \frac{\eta(3\t)^{9}}{\eta(\t)^3}=\sum_{n=1}^\inf
\bigl(\sum_{d|n} \left(\frac{d}{3}\right)(n/d)^2\bigr) q^n $$
have
non-negative and positive Fourier coefficients respectively, as shown 
below. Notice that the coefficients are
multiplicative and hence checking when $n$ is a power of prime is
enough. If $n=p^e$ is a power of prime, we have
$$1+\left(\frac{p}3\right)+\left(\frac{p^2}3\right)+\cdots
+\left(\frac{p^e}3\right)
=\cases 1, &\qquad p=3,\\
e+1, &\qquad \left(\dfrac{p}3\right)=1,\\
1, &\qquad \left(\dfrac{p}3\right)=-1\ \text{and}\ e=\text{even},\\
0, &\qquad \left(\dfrac{p}3\right)=-1\ \text{and}\ e=\text{odd}, \endcases
$$ and
$$p^{2e}+\left(\frac{p}3\right)p^{2e-2}
+\left(\frac{p^2}3\right)p^{2e-4}+\cdots
+\left(\frac{p^e}3\right)
=\cases p^{2e}, &\quad p=3,\\
\dfrac{p^{2e+2}-1}{p^2-1}, &\quad \left(\dfrac{p}3\right)=1,\\
\dfrac{p^{2e+2}+1}{p^2+1}, &\quad \left(\dfrac{p}3\right)=-1\ \text{and}
\ e=\text{even},\\
\dfrac{p^{2e+2}-1}{p^2+1}, &\quad \left(\dfrac{p}3\right)=-1\ \text{and}
\ e=\text{odd}. \endcases $$ 
From this we conclude that $\sum_{d|n} \left(\dfrac{d}{3}\right)\ge 0$
and $\sum_{d|n} \left(\dfrac{d}{3}\right)(n/d)^2>0$.

When $k$ is of the form $(6n+1)/2$, the product
$(-(2k-1)/6)_i(-(k-2)/6)_i$ needs not always be positive and so we
need some extra work.  Put $$a_i=
\frac{(-\frac{2k-1}{6})_i(-\frac{k-2}{6})_i} {(\frac{k+7}6)_i
  \,i!}=\frac{(-n)_i(-\frac{2n-1}4)_i}{(\frac{2n+5}{4})_i \,i!}.$$
We
have to show that the form
$$\D_2^{(4)}(\t)^{\frac{2n+1}4}\sum_{i=0}^n a_i
(16\D_2^{(4)}(\t))^{i}E_2^{(4)}(\t)^{n-i}$$
has  positive Fourier
coefficients. Put $A=\D_2^{(4)}(\t)$ and $B=E_2^{(4)}(\t)$ as in the
proof of Theorem 1. Note that the number $a_i$ is positive for
$i<(2n+3)/4$ and beyond this the sign alternates. Suppose now
$a_i$ is positive and  $a_{i+1}$ negative (or zero). Since
$\D_2^{(4)}=\sum_{n\ge 1,\text{odd}} \bigl(\sum_{d|n}d\bigr) q^n$ has
positive coefficients,  our proof is complete if we show that the
Fourier coefficients of
$$\align &a_i (16A)^{i}B^{n-i}+
a_{i+1} (16A)^{i+1}B^{n-i-1}\\&=a_i(16A)^{i}B^{n-i-1}\left(B+
16\frac{(-n+i)(-\frac{2n-1}4+i)}{(\frac{2n+5}4+i)(i+1)}A\right)\endalign $$ 
are positive. For this we prove the following lemma.

\proclaim{Lemma} For $\a$ with $0\le \a <8$, the Fourier coefficients
in $E_2^{(4)}(\t)-\a \D_2^{(4)}(\t)$ are positive. \endproclaim

\pf Denote the sum $\sum_{d|n}d$ by $\s(n)$.  By
$$E_2^{(4)}(\t)=\frac13\left(4E_2(4\t)-E_2(\t)\right)
=1+8\sum_{n=1}^\inf \s(n)q^n-32\sum_{n=1}^\inf \s(n)q^{4n}$$
and
$$\D_2^{(4)}(\t)=\sum_{n\ge 1,\text{odd}}\s(n)q^n,$$
we have
$$\align &E_2^{(4)}(\t)-\a \D_2^{(4)}(\t)\\
&=1+8\sum_{n\equiv 2 \bmod4}\s(n)q^n+
8\sum_{n=1}^\inf\left(\s(4n)-4\s(n)\right)q^{4n}
+\sum_{n\ge 1,\text{odd}}(8-\a)\s(n)q^n.
\endalign$$  To see the coefficient in the middle sum is positive,
write $n=2^em$, $m$:\,odd. Then we have
$$\align \s(4n)-4\s(n)&=\s(2^{e+2}m)-4\s(2^em)=(\s(2^{e+2})-4\s(2^e))\s(m)\\
&=\left(2^{e+3}-1-4(2^{e+1}-1)\right)\s(m)=3\s(m)>0. \endalign$$ 
Thus we have proved the lemma. \endpf

In view of the lemma, the proof of Theorem 3 is complete if the inequality
$$0\le -16\frac{(-n+i)(-\frac{2n-1}4+i)}{(\frac{2n+5}4+i)(i+1)}<8\quad
\text{for}\quad \frac{2n+3}4<i\le n$$
is shown  to hold, which is
readily seen. \endpf

\section{\S5. A characterization of the equation \eq}

Consider the differential equation in the upper half-plane of the form
$$
f''(\t)+A(\t)f'(\t)+B(\t)f(\t)=0 \tag1$$
where $A(\t)$ and
$B(\tau)$ are assumed to be holomorphic in $\h$ and  bounded when
$\Im(\t)\rightarrow\inf$. Fix a non-negative integer $k$ and we
further assume: \medskip

If $f(\t)$ is a solution of (1), then $(c\t+d)^{-k}f\left(\dfrac
  {a\t+b}{c\t+d}\right)$ is also a solution  of (1)

for all $\pmatrix
a & b \\ c & d \endpmatrix \in\g$.
\medskip

\proclaim{Proposition 2} The differential equation \eq is essentially
the unique equation which satisfies the above conditions. \endproclaim

\pf Put $g(\t)=(c\t+d)^{-k}f\left(\dfrac{a\t+b}{c\t+d}\right)$. By an
elementary calculation, we have from the assumption
$$ \align 0&= (c\t+d)^{k+4}\left(g''(\t)+A(\t)g'(\t)+B(\t)g(\t)\right)\\
&= f''\left(\frac{a\t+b}{c\t+d}\right)+\left((c\t+d)^2A(\t)-
\frac{k+1}{\pi i}c(c\t+d)\right)f'\left(\frac{a\t+b}{c\t+d}\right)\\
&+\left((c\t+d)^4B(\t)-\frac{k}{2\pi i}c(c\t+d)^3A(\t)+\frac
{k(k+1)}{(2\pi i)^2}c^2(c\t+d)^2\right)f\left(\frac{a\t+b}{c\t+d}\right).
\endalign $$ Comparing this with 
$$f''\left(\frac{a\t+b}{c\t+d}\right)+A\left(\frac{a\t+b}{c\t+d}\right)
f'\left(\frac{a\t+b}{c\t+d}\right)
+B\left(\frac{a\t+b}{c\t+d}\right)f\left(\frac{a\t+b}{c\t+d}\right)=0,$$
we have (under the natural assumption that $f$ and $f'$ are
independent)
$$A\left(\frac{a\t+b}{c\t+d}\right)=(c\t+d)^2A(\t)-
\frac{k+1}{\pi i}c(c\t+d) \tag2$$ and 
$$B\left(\frac{a\t+b}{c\t+d}\right)=(c\t+d)^4B(\t)-\frac{k}{2\pi i}
c(c\t+d)^3A(\t)+\frac{k(k+1)}{(2\pi i)^2}c^2(c\t+d)^2.\tag3$$
By the
quasimodular property of $E_2(\t)$ (recalled in Section 2), the equation (2)
says that $A(\t)+(k+1)E_2(\t)/6$ transforms like a modular form of
weight $2$ under the action of $\g$. By our assumption that $A(\t)$ is
holomorphic on $\h$ and bounded when $\Im(\t)\rightarrow\inf$,  this
must be $0$ and hence
$$A(\t)=-\frac{k+1}6E_2(\t).$$ Putting this into (2) and using 
$$E_2'\left(\frac{a\t+b}{c\t+d}\right)=(c\t+d)^4E_2'(\t)+
\frac1{\pi i}c(c\t+d)^3E_2(\t)+\frac3{(\pi i)^2}c^2(c\t+d)^2,$$
we conclude that $B(\t)-k(k+1)E_2'(\t)/12$ behaves like
a modular form of weight $4$ on $\g$. Hence, we have 
$$B(\t)=\frac{k(k+1)}{12}E_2'(\t) +\alpha E_4(\t)$$
with some constant
$\a$. But an easy calculation shows that if  $f(\t)$ is a solution of
$$f''(\t)-\frac{k+1}6E_2(\t)f'(\t)
+\left(\frac{k(k+1)}{12}E_2'(\t)+\alpha E_4(\t)\right)f(\t)=0,\tag5$$
then $f(\t)\D(\t)^\b$ is a solution of $(\#)_{k'}$ with $\b$ a
solution of $$\b^2+\frac{k+1}6\b+\a=0$$
and $k'=k+12\b.$
Hence we may take $\a=0$ without losing any generality and we therefore
conclude that \eq is  essentially unique.  \endpf

\example{Remark} The above quadratic equation for $\b$ has no real
solution when $$\left(\frac{k+1}6\right)^2-4\a<0.\tag6$$
However, if
the equation (5) has a power series solution starting with
$q^m+\cdots$, then by equating the coefficient of $q^m$ we should have
the relation
$$m^2-\frac{k+1}6m+\a=0.$$
Thus, if the inequality (6) holds, no power series solution to (5) can exist. 
\endexample \medskip

If we loosen the condition posed on the coefficients of (1), we get a
variety of similar differential equations. It may be an interesting
task to seek modular or quasimodular solutions of those equations
and study their properties.  \medskip

The $\g$-invariance property of solutions of \eq has an interesting
consequence concerning the space of solutions when $k\equiv 5\mod 6$
(quasimodular case) and the representation of $\g$ on this space.

We first observe that the quasimodular solution in Theorem 2 can be
written in the form $$f(\t)E_2(\t)+g(\t)\ (=:F_k(\t))$$
where $f(\t)$
and $g(\t)$ are $\g$-modular of weight $k-1$ and $k+1$
respectively. (We have used $E_4'=(E_4E_2-E_6)/3.$) Write the action
of $\g$ of weight $k$ on a  function $F$ by
$$F \vert_k \gamma := (c\t+d)^{-k}F\left(\frac{a\t+b}{c\t+d}\right)
\quad\text{for}\quad \gamma=\pmatrix a & b \\ c & d \endpmatrix \in\g.$$
By the transformation formula 
$$E_2\left(\frac{a\t+b}{c\t+d}\right)=(c\t+d)^2E_2(\t)+\frac6{\pi i}
c(c\t+d),$$
we have $$
\align & F_k(\t)\vert_k \left(\smallmatrix a & b \\
  c & d \endsmallmatrix\right)\\
&=(c\t+b)^{-k}\biggl((c\t+b)^{k-1}f(\t)
\bigl((c\t+b)^2E_2(\t)+\frac6{\pi i}c(c\t+b)\bigr)+(c\t+b)^{k+1}g(\t)\biggr)\\
&=(c\t+b)\bigl(f(\t)E_2(\t)+g(\t)\bigr)+\frac6{\pi i}cf(\t)\\
&=c\bigl(\t F_k(\t)+\frac6{\pi i}f(\t)\bigr)+dF_k(\t).\tag7 \endalign
$$
Hence we see that the function
 $$\t F_k(\t)+\frac6{\pi i}f(\t)$$ is also a solution of $(\#)_k$ and
the space of solutions is spanned by this and the $F_k(\t)$. Moreover,
the action of $\g$ on this function is computed as
$$\align &\bigl(\t F_k(\t)+\frac6{\pi i}f(\t)\bigr)\vert_k
\left(\smallmatrix a & b \\
  c & d \endsmallmatrix\right)\\
&=\frac{a\t+b}{c\t+d}\bigl(F_k(\t)\vert_k \left(\smallmatrix a & b \\
  c & d \endsmallmatrix\right)\bigr)+\frac6{\pi i}(c\t+d)\i f(\t)\\
&=\frac{a\t+b}{c\t+d}\bigl((c\t+d)F_k(\t)+\frac6{\pi i}cf(\t)\bigr)
+\frac6{\pi i}(c\t+d)\i f(\t)\\
&=(a\t+b)F_k(\t)+\frac6{\pi i}af(\t)\\
&=a\bigl(\t F_k(\t)+\frac6{\pi i}f(\t)\bigr)+bF_k(\t).\endalign $$
This, together with (7), shows that the representation of $\g$ on the
space of solutions of \eq in the case $k\equiv 5\mod 6$ is  faithful
and in particular identical if we choose  $\t F_k(\t)+6(\pi i)\i
f(\t)$ and $F_k(\t)$ as a basis. 

As a corollary, we conclude that when $k\equiv 5\bmod 6$, no modular
solution of \eq of weight $k$ on any subgroup of $\g$ exists because
otherwise there would exist a subspace invariant under the action of
that group. As for the  other cases, we have no proof that the
solutions presented in Theorem 1  exhaust all possible modular
solutions (at least when $k$ is an integer or half an integer) but
conjecture so on the basis of numerical  experiments. As mentioned in the
end of Section 3, it would be desirable to have, if any, a number
theoretical reason for the existence of modular solutions.

\section{BIBLIOGRAPHY}

\ref\no 1  \by  M.~Kaneko and D.~Zagier \pages 97--126  \paper
Supersingular $j$-invariants,  Hypergeometric series, and Atkin's
orthogonal polynomials \yr 1998\vol 7\jour AMS/IP Studies in Advanced
Mathematics \endref

\ref\no 2  \by  Ikuo Satake \pages 247--251 \paper
Flat structure for the simple elliptic singularity of type
$\tilde{E}_6$ and Jacobi form\yr 1993 \vol 69, Ser. A, No. 7
\jour Proc. of the Japan Academy
\endref

\ref\no 3  \by  Ikuo Satake \pages 427--452 \paper  Flat structure and the 
prepotential for the elliptic root system of type $D^{(1,1)}_4$, {\rm in
Topological Field Theory, Primitive Forms and Related Topics
(Kashiwara, Matsuo, Saito, and Satake eds.)} 
\yr 1998 \vol 160 \jour Progress in Math.
\endref

\end